\documentclass{article}
\usepackage{amsmath,amsfonts,amssymb,amsthm,amscd}
\title{O-minimal spectra, infinitesimal
subgroups and cohomology} 
\author{Alessandro Berarducci\footnote{URL:
www.dm.unipi.it/$\sim$berardu. Partially supported by Progetto MIUR,
Cofin 2004, Metodi di Logica in Algebra, Analisi e Geometria. Part of
the results were presented at the meeting ``Around o-minimality'',
March 11-13, 2006, Leeds.} \\Universit\`a di Pisa\\Largo Bruno Pontecorvo
5\\56127 Pisa}
\date{April 9, 2006}  \parindent=0mm
\DeclareMathOperator{\R}{\mathbb R}
\DeclareMathOperator{\NN}{\mathbb N}

\DeclareMathOperator{\Z}{\mathbb Z}

 \DeclareMathOperator{\sta}{st}
\DeclareMathOperator{\df}{df} \DeclareMathOperator{\stab}{Stab}
\DeclareMathOperator{\Def}{Def} 
 
\DeclareMathOperator{\Cl}{Cl}
\DeclareMathOperator{\Int}{Int}
\newcommand{\ov}{\overline}
\newcommand{\nin}{\not\in}

\newcommand{\nog}{{\cal I}}

\theoremstyle{plain}
\newtheorem{theorem}{Theorem}
\newtheorem{lemma}[theorem]{Lemma}
\newtheorem{proposition}[theorem]{Proposition}
\newtheorem{corollary}[theorem]{Corollary}
\newtheorem{conjecture}[theorem]{Conjecture}
\newtheorem{claim}{Claim}

\theoremstyle{definition}
\newtheorem{remark}[theorem]{Remark}
\newtheorem{definition}[theorem]{Definition}
\newtheorem{example}[theorem]{Example}

\numberwithin{theorem}{section}
\begin{document}
\maketitle

\begin{abstract} 
By recent work on some conjectures of Pillay, each definably compact
group $G$ in a saturated o-minimal expansion of an ordered field has
a normal ``infinitesimal subgroup'' $G^{00}$ such that the quotient
$G/G^{00}$, equipped with the ``logic topology'', is a compact (real)
Lie group.  Our first result is that the functor $G\mapsto G/G^{00}$
sends exact sequences of definably compact groups into exacts
sequences of Lie groups. We then study the connections between the Lie
group $G/G^{00}$ and the o-minimal spectrum $\widetilde G$ of $G$. We
prove that $G/G^{00}$ is a topological quotient of $\widetilde G$. We
thus obtain a natural homomorphism $\Psi^*$ from the cohomology of
$G/G^{00}$ to the (\v{C}ech-)cohomology of $\widetilde G$. We show
that if $G^{00}$ satisfies a suitable contractibility conjecture then
$\widetilde {G^{00}}$ is acyclic in \v{C}ech cohomology and $\Psi^*$
is an isomorphism. Finally we prove the conjecture in some special cases.
\end{abstract} 

\section{Introduction}
We discuss some topics related with Pillay's conjectures
\cite{Pillay04} for groups definable in an o-minimal expansions of a
field.  We begin with an example which illustrates the situation in a
simplified setting.

\begin{example} \label{Example} Let $SO(n,\R)$ be the special orthogonal
group, a compact connected Lie group. Let $M$ be a real closed
field properly containing $\R$ and consider the group $SO(n,M)$.
Each element of $SO(n,M)$ has coordinates bounded by an element of
$\R$ (indeed they are all $\leq 1$ in absolute value).  Thus there
is a surjective group homomorphism $\sta \colon SO(n,M) \to
SO(n,\R)$ where $\sta(x)$ is the unique $y\in SO(n,\R)$ whose
coordinates satisfy $|x_i - y_i| < 1/k$ for each positive integer
$k$.

Let $SO(n,M)^{00}$ be the kernel of $\sta \colon SO(n,M) \to
SO(n,\R)$. We call $SO(n,M)^{00}$ the ``infinitesimal subgroup'' of
$SO(n,M)$. The quotient $SO(n,M)/SO(n,M)^{00}$ is isomorphic to
$SO(n,\R)$ as an abstract group, and also as a topological group
provided we give it the following ``logic topology'' (which is not a
quotient topology): a subset $X \subset SO(n,M)/SO(n,M)^{00}$ is
closed iff its preimage in $SO(n,M)$ is the intersection of a
countable family of semialgebraic sets.
\end{example}

The following provides an intrinsic characterization of $SO(n,M)^{00}$
which is invariant under semialgebraic group isomorphisms over $M$ (unlike the
one using the map $\sta$).

\begin{proposition} \label{fact} $SO(n,M)^{00}$ is the smallest subgroup of $SO(n,M)$
which has index $\leq 2^{\aleph_0}$ and is the intersection of a
countable family of semialgebraic sets. \end{proposition}

\begin{proof} By Definition \ref{G00}, Theorem \ref{good} and Lemma
\ref{ctble}. \end{proof}

It turns out that it is possible to introduce a notion of
infinitesimal subgroup and generalize the above observations in a more
general context.  We consider definable groups $G\subset M^n$ in a big
saturated o-minimal expansion $M$ of a real closed field. For a recent
survey on definable groups see \cite{Otero06}. The following results give a
positive answer to some conjectures of Pillay in \cite{Pillay04}. For
the definitions of the notions involved (logic-topology,
type-definable, bounded index, definably compact) see section 2.

\begin{theorem} (\cite{BerarducciOPP05}) \label{BOPP04} Let $G$ be a definable
group. Then $G$ has the descending chain condition on type-definable
subgroups of bounded index. If $H\lhd G$ is a normal type-definable
subgroups of bounded index, then $G/H$ with the logic topology is a
compact (real) Lie group. \end{theorem}

\begin{definition} \label{G00} By the descending chain condition (or by
\cite{Shelah05}) there exists a smallest type-definable subgroups
$G^{00}<G$ of bounded index, which is in fact necessarily normal in
$G$ \cite[Remark 2.9]{Pillay04}. \end{definition}

Under the assumption that $G$ is definably compact in
the sense of \cite{PeterzilS99} we obtain.

\begin{theorem} (\cite{HrushowskiPP05}) \label{HPP05} Let $G$ be a definably
compact definable group.  Then the dimension of $G/G^{00}$ as a Lie
group equals the o-minimal dimension of $G$.  \end{theorem}

The following result handles the situation that we encountered in
Example \ref{Example}. 

\begin{theorem} (\cite[Prop. 3.6]{Pillay04} and \cite[Fact 4.1]{PeterzilP04})
\label{good} Let $G$ be a definably compact definable group.  If $M$
contains as an elementary substructure an o-minimal expansion of the
real field $\R$, and if $G$ is defined over $\R$, then $G/G^{00}
\simeq G(\R)$ as Lie groups. \end{theorem}

The situation of Theorem \ref{good} is not the typical one. In fact, even assuming
that $M$ contains $\R$ as an elementary substructure, there are
definable groups over $M$ which are not definably isomorphic to groups
defined over $\R$ (this happens even for elliptic curves
\cite{PeterzilS04}), so in general we cannot naturally identify
$G/G^{00}$ with some group defined over $\R$. Nevertheless $G/G^{00}$
is always a compact Lie group.

Let us now discuss the contributions of this note.  A point of view
which until now has not been stressed is that the correspondence $G
\mapsto G/G^{00}$ is a functor from definably compact groups and
definable group homomorphisms to compact Lie groups and continuous
homomorphism.  After giving an exposition of the relevant results from
\cite{PeterzilP04,HrushowskiPP05}, we show that this functor is exact,
namely it sends short exact sequences $0\to H\to G\to B \to 0$ to
short exact sequences $0\to H/H^{00} \to G/G^{00} \to B/B^{00} \to
0$. Many of the arguments are implicit in \cite{HrushowskiPP05}. The
extra ingredient needed is a conjugacy theorem for maximal tori proved
in \cite{Edmundo05} (see also \cite{Berarducci05} for a different
proof).  As a side-product we characterize $G^{00}$ as the unique
type-definable subgroup of $G$ of bounded index which is divisible and
torsion free (this was proved in \cite{HrushowskiPP05} if $G$ is
abelian and definably connected).

The exactness of the functor $G\mapsto G/G^{00}$ may find applications
in the attempt to compare the topological invariants of $G$ and
$G/G^{00}$. This idea however is still to be explored.

In the second part of the paper we study the connection between
$G/G^{00}$ and the space $\widetilde G$ consisting of the types of $G$
(over a small model $M_0\prec M$) with the spectral topology. The
spectral space $\widetilde G$ is an o-minimal analogue of the real
spectrum of real semialgebraic geometry
\cite{CosteR82,CarralC83,Pillay88} and it is a quasi-compact normal
topological space. Note that $\widetilde G$ is not Hausdorff and does
not carry a group structure. There is a natural surjective map $\Psi
\colon \widetilde G \to G/G^{00}$ sending a type $p$ to the coset
$gG^{00} \in G/G^{00}$ of a realization $g$ of $p$. We show that
$\Psi$ is a continuous closed map, and therefore it is a quotient map,
namely a subset of $U$ of $G/G^{00}$ is open in the logic topology if
and only if $\Psi^{-1}(U)$ is open in the spectral topology. Thus we
can identify the topological space $G/G^{00}$ (with the logic
topology) as a quotient of $\widetilde G$ modulo the equivalence
relation $\ker (\Psi)$ (with the quotient topology).

The map $\Psi$ induces an homomorphism $\Psi^*\colon
\check{H}^n(G/G^{00}) \to \check{H}^n(\widetilde G)$ in \v{C}ech
cohomology. In the last part of the paper we show that if $G^{00}$
satisfies a certain contractibility property, then $\widetilde
{G^{00}}$ is acyclic in \v{C}ech cohomology and $\Psi^*$ is an
isomorphism.  In this case we also obtain a natural isomorphism
${H}^n(G/G^{00}) \simeq {H}_{\df}^n( G)$ where $H^n$ denotes
the singular cohomology functor and $H^n_{\df}$ is the definable
singular cohomology functor studied in \cite{EdmundoO04} and based on
the definable singular homology of \cite{Woerheide96}.

\section{Logic topology}
As in \cite{Pillay04,BerarducciOPP05}, we adopt the model theoretic
convention of working in a ``universal domain''. So let $M$ be a
saturated o-minimal $L$-structure of cardinality $\kappa$, where
$\kappa$ is say inaccessible, and strictly greater than the
cardinality of the language of $M$. By a {\bf definable set} in $M$ we
mean a subset of some $M^{n}$ defined by an $L$-formula with
additional parameters from $M$. By a {\bf type-definable set} in $M$
we mean a subset of $M^{n}$ which is the intersection of
$<\kappa$-many definable sets. If $A\subset M$, we say that $X$ is
(type-)definable over $A$ if $X$ can be defined by a formula (set of
formulas) with parameters from $A$. If $X$ is a (possibly
type)-definable set and $E$ is a type-definable equivalence relation
on $X$, then we call $E$ {\bf bounded} (or of {\bf bounded index}) if
$|X/E| < \kappa$. The real meaning of boundedness is that $X/E$ does
not increase when we pass to an elementary extension $M'\succ M$. More
precisely the canonical injection $i \colon X/E \to X(M')/E(M')$,
$i(a/E) = a/E(M')$, is a bijection \cite[Fact 2.1]{Pillay04}. If $E$
is a bounded type-definable equivalence relation on $X$, we put on
$X/E$ the {\bf logic topology}, in which a subset $Z\subseteq X/E$ is
defined to be closed if its preimage under the natural map $\pi\colon
X\to X/E$ is type-definable.  Then $X/E$ is a compact (Hausdorff)
topological space \cite[Lemma 2.5]{Pillay04}.

In the sequel by a {\bf small} set we mean a set of cardinality
$<\kappa$. So if $M_0\prec M$ is a small model, each type over $M_0$
is realized in $M$. 

\begin{lemma} \label{repres} (\cite[Remark 1.6]{BerarducciOPP05}) Let $M_0\prec
M$ be a small model over which $X,E$ are (type)-defined. Then
$Z\subset X/E$ is closed iff its preimage under $\pi\colon X \to X/E$
is type-definable over $M_0$. \end{lemma}

The proof of Lemma \ref{repres} depends on \cite{LascarP01} but there
is an easier proof if we further assume, as we will do in the sequel,
that $M_0$ contains a representative from each equivalence class (this
can be done keeping $M_0$ small since $E$ is bounded).  In fact we can
put on $X/E$ two topologies: the logic topology, and the one which
declares a set closed if its preimage in $X$ is type-definable over
$M_0$. The proof in \cite[Lemma 2.5]{Pillay04} that $X/E$ is compact
Hausdorff applies to both topologies (the extra assumption is used to
prove the Hausdorff property), and since one is finer than the other
they must coincide.

We will be mainly interested in the case when $X$ is a definable group
$G$ and $E$ is the equivalence relation induced by the left-cosets of
a type-definable subgroup $H<G$. We say that $H$ is bounded if $H$ has
a small number of left cosets in $G$. In this case we can put on
$G/H$ the logic topology.  So in Example \ref{Example} the subgroup
$SO(n,M)^{00}$ of $SO(n,M)$ is bounded since it has index
$2^{\aleph_0}$ ($=$ the cardinality of $SO(n,M)/SO(n,M)^{00} \simeq
SO(n,\R)$).  In general by Theorem \ref{BOPP04} every normal bounded
type-definable subgroup $H$ of a definable group $G$ has index $\leq
2^{\aleph_0}$ since the quotient $G/H$ is a Lie group. (The same
conclusion holds if $H$ is not normal since every type-definable
bounded subgroup contains a normal type-definable subgroup of bounded
index \cite{Pillay04}.)

By \cite{Pillay88b} any definable group $G$ has a unique topology
which makes it into a ``definable manifold''. Since we are working
over an o-minimal expansion of a field, by the o-minimal version of
Robson's embedding theorem \cite[Ch. 10, Thm. 1.8]{vdDries98}, $G$ is
definably homeomorphic to a definable submanifold of $M^n$, where $M$
has the topology generated by the open intervals and $M^n$ the product
topology. (The embedding theorem works under very general hypothesis,
but for the case of definably compact groups there are simpler proofs,
see \cite[Thm. 10.7]{BerarducciO01}.) So we can assume without loss of
generality that $G$ is a submanifold of $M^n$ and the group operation
is continuous in the induced topology from $M^n$.  Under this
assumption $G$ is {\bf definably compact} in the sense of
\cite{PeterzilS99} iff $G$ is closed and bounded in $M^n$. A definable
set $X\subset M^n$ is {\bf definably connected} if $X$ cannot be
written as the disjoint union of two definable non-empty open
subsets. If $G$ is a definable group (embedded in $M^n$), $G$ is
definably connected iff it has no proper subgroups of finite index
\cite[Prop. 2.12]{Pillay88b}.  In general if $X$ is a subset of $M^n$
we put on $X$ the topology induced by $M^n$.

\begin{lemma} \label{ctble} Let $G$ be a definable group and let $H\lhd G$ be a
type-definable normal subgroup of $G$ of bounded index. Then $H$ is
the intersection of a countable decreasing sequence $Y_i$ of definable
subsets of $G$. Moreover we can arrange so that $\Cl(Y_{i+1})\subset
\Int(Y_i)$, where $\Cl$ and $\Int$ denote the closure and the interior
in $G$.  \end{lemma}

\begin{proof} Being a Lie group $G/H$ admits a decreasing sequence $O_1 \supset
O_2 \supset O_3 \supset \ldots$ of open neighbourhoods of the neutral
element $e$ with $\bigcap_i O_i = \{e\}$ and $\ov{O_{i+1}} \subset
O_i$. The preimages under $\pi \colon G \to G/G^{00}$ of the closed
sets $\ov{O_{i+1}}$ and the open sets $O_i$ are respectively
type-definable and $\bigvee$-definable (small unions of definable
sets). In a saturated model if a type-definable set is contained in a
$\bigvee$-definable set, there is a definable set in between. So there
are definable sets $X_i$ with $\pi^{-1}(\ov{O_{i+1}}) \subset X_i
\subset \pi^{-1}(O_{i})$ and we obtain $\bigcap_i X_i = H$. Since
$\pi$ is continuous (Lemma \ref{cont}) $\pi^{-1}(\ov{O_{i+1}})$ is a
closed subset of the open set $\pi^{-1}(O_{i})$ and
$\Cl(X_{i+1})\subset \pi^{-1}(\ov{O_{i+1}}) \subset \pi^{-1}(O_i)
\subset \Int(X_{i-1})$, so we can set $Y_i = X_{2i}$. \end{proof}

\begin{remark} If $G$ is defined over $M_0\prec M$, then $G^{00}$ is
type-definable over $M_0$ (\cite{Shelah05} or
\cite[Prop. 6.1]{HrushowskiPP05}) and by Lemma \ref{repres} we can
take the sets $X_i$ to be $M_0$-definable. \end{remark}

\section{Generic sets}

In this section we recall some results and definitions from
\cite{PeterzilP04,HrushowskiPP05}.  Let $G$ be a definably compact
group (in $M$) and let $\Def(G)$ be the family of all definable
subsets of $G$.  A set $X\in \Def(G)$ is {\bf generic} iff finitely
many left-translates of $X$ cover $G$. This is equivalent to require
that finitely many right translates of $X$ cover $G$
\cite[Prop. 4.2]{HrushowskiPP05}. If $X\in \Def(G)$ contains a
type-definable subgroup of bounded index, then it is easy to see using
the saturation of $M$ that $X$ is generic. So in particular if
$X\supset G^{00}$, then $X$ is generic. If the union of two definable
subsets of $G$ is generic, one of the two is generic \cite[Prop. 4.2,
Thm. 8.1]{HrushowskiPP05}. So the non-generic definable subsets of $G$
form an ideal ${\cal I}$. Consequently we can introduce an equivalence
relation $\sim_{\nog}$ on $\Def(G)$ by setting $X \sim_{\nog} Y$ iff
the symmetric difference $X\triangle Y$ is not generic.  Clearly if $g
\in G$ and $X\sim_{\nog} Y$ then $gX \sim_{\nog} gY$. Therefore for
each $X\in \Def(G)$ the set
\[\stab_{\nog}(X) = \{g\in G \mid gX \sim_{\nog} X\}\eqno (1) \] is a subgroup
of $G$.  It is immediate to see that $\stab_{\nog}(X)$ is
type-definable. One of the main results of \cite{HrushowskiPP05} is:

\begin{theorem} (\cite{HrushowskiPP05}) \label{bounded} $\stab_{\nog}(X)$ has
bounded index in $G$. \end{theorem}

This is a consequence of the following: 

\begin{lemma} (\cite[Remark. 4.4, Thm. 8.1]{HrushowskiPP05}) If $M_0\prec M$ is
a model over which $G$ is defined, then every generic set $X\subset G$
(not necessarily defined over $M_0$) meets $M_0$. \end{lemma} 

Most of the results in \cite{HrushowskiPP05} are stated for simplicity
under the assumption that $G$ is definably compact and definably
connected, but it is easy to see that in the above results it suffices
that $G$ is definably compact.  Granted the lemma, to prove Theorem
\ref{bounded} one fixes a small model $M_0\prec M$ over which $G$ is
defined. If $g,h\in G$ have the same type over $M_0$, then the set $gX
\triangle hX$ does not meet $M_0$, and therefore $gX \sim_{\nog}
hX$. Now $g\stab_{\nog}(X) = h\stab_{\nog}(X)$ iff $gX \sim_{\nog}
hX$.  So the index of $\stab_{\nog}(X)$ is bounded by the cardinality
of $S_G(M_0)$, where $S_G(M_0)$ is the set of types over $M_0$ of
elements of $G$.

An immediate consequence of 
Theorem \ref{bounded} is:

\begin{corollary} (\cite{HrushowskiPP05}) \label{stab} $G^{00} = \bigcap_{X\in
\Def(G)} \stab_{\nog} (X)$. \end{corollary}

\begin{proof} The $\subset$ inclusion follows from Theorem \ref{bounded}.  For
the opposite inclusion suppose $gX \sim_{\nog} X$ for all $X\in
\Def(G)$. We want to show $g\in G^{00}$. If not $gG^{00} \cap G^{00} =
\emptyset$. Since $G^{00}$ is a small intersection of definable sets,
by saturation there is a definable set $X\supset G^{00}$ such that $gX
\cap X = \emptyset$. Together with $gX \sim_{\nog} X$ this implies
$X\in \nog$. This is absurd since $X \supset G^{00}$. \end{proof}

We will later need the following: 

\begin{corollary} \label{CC} Let $G$ be a definably compact group. Then \[G^{00} =
\bigcap_{X\mbox{ generic }}XX^{-1}.\] \end{corollary}

\begin{proof} Let $X\in \Def(G)$ be generic. Then $G^{00} \subset
\stab_{\nog}(X) \subset \{g \in G \mid gX \cap X \neq \emptyset\} =
XX^{-1}$. So $G^{00} \subset \bigcap_{X\in \Def(G)}XX^{-1}$. For the
other inclusion write $G^{00} = \bigcap_i X_i$ where $\{X_i \mid i\in
I\}$ is a downward directed small family of definable sets and note
that by saturation $(\bigcap_i X_i)(\bigcap_i X_i)^{-1} = \bigcap_i X_iX_i^{-1}$. \end{proof}

\begin{remark} If $G$ is definable over $M_0$, in Corollary \ref{stab} and
Corollary \ref{CC} we can restrict $X$ to range over the $M_0$-definable sets. \end{remark}

The sets $X\in \Def(G)$ whose complement is non-generic are stable
under finite intersections.  It follows that there is a type $p\in
S_G(M)$ containing all the sets $X\in \Def(G)$ whose complement is
non-generic. Clearly for all $X\in \Def(G)$, if $X\in p$ then $X$ is
generic. We call a type with this property a {\bf generic type}. So
each definably compact group $G$ has a generic type. It is shown in
\cite{HrushowskiPP05} that each generic type $p \in S_G(M)$ induces a
left-invariant finitely additive measure on $\Def(G)$ by
\[\mu_p(X) = m\{gG^{00} \in G/G^{00} \mid gX \in p\} \eqno (2) \] where $m$ is the Haar measure on the compact Lie group $G/G^{00}$. 
The set $\{gG^{00} \mid gX \in p\}$ is well defined by Corollary
\ref{stab}. Indeed if $gG^{00} = hG^{00}$, then $gh^{-1}\in G^{00}
\subset \stab_{\nog}(X)$, so $gX \sim_{\nog} hX$ and therefore $gX\in
p \leftrightarrow hX\in p$. Moreover by \cite{HrushowskiPP05}
$\{gG^{00} \in G/G^{00} \mid gX \in p\}$ is a Borel set, so it has a
well defined Haar measure. Thus we have:

\begin{theorem} (\cite{HrushowskiPP05}) There is a finitely additive left
invariant measure \[\mu \colon \Def(G) \to [0,1]^{\R} \] with $\mu(G)
= 1$. Moreover for $X\in \Def(G)$, $\mu(X)>0$ iff $X$ is generic. \end{theorem}

For the last statement note that if $X$ is not generic, then
$\{gG^{00} \mid gX \in p\}$ is empty, and therefore $\mu(X) = 0$.

\begin{conjecture} 
There is a unique finitely additive measure $\mu \colon \Def(G)\to [0,1]$ with $\mu(G) = 1$. 
\end{conjecture} 

Indeed one may even conjecture that $\mu(X) = m(\pi(X))$ where
$\pi(X)$ is the image of $X\subset G$ in $G/G^{00}$. This would imply
in particular that if $X$ has empty interior in $G$, then $m(\pi(X)) =
0$. By \cite[Lemma 10.5]{HrushowskiPP05} this in turn implies the
compact domination conjecture of that paper. A possible attempt to
prove the conjecture could be based on the ideas of \cite[Thm. A, \S 13, p.
54]{Halmos74} where it is proved that a ``measurable group'' has a
unique measure up to a constant factor.

\section{Torsion free subgroups} 

\begin{lemma} \label{torsionfreeabelian} Let
$G$ be a definably compact definably connected abelian group. Then
$G^{00}$ is the unique normal type-definable subgroup of bounded index
which is torsion free. \end{lemma}

\begin{proof} $G^{00}$ is normal and divisible by \cite[Fact 1.8 and \S 3, claim
3]{BerarducciOPP05}. By \cite[Proof of 3.15 ]{PeterzilP04} for each
$n$ one can find a generic $X\subset G$ such that $\stab_{\nog}(X)$
has no $n$-torsion. By \cite{HrushowskiPP05} $G^{00} = \bigcap_X
\stab_{\nog}(X)$, so $G^{00}$ is torsion free. Conversely if $H$ is a
torsion free type-definable subgroup of $G$ of bounded index, then
$H=G^{00}$ by \cite[Cor. 1.2]{BerarducciOPP05}. \end{proof}

We next investigate the effect of dropping the assumption that $G$ is
abelian and definably connected. Recall that every definable group $G$
has a smallest definable subgroup $G^0$ of finite index
\cite{Pillay88b}. Then clearly $G^0$ is the unique definably connected
definable subgroup of $G$. Moreover $G^{00} < G^0$ and we easily obtain:

\begin{lemma} \label{G0} $(G^0)^{00} = G^{00}$. \end{lemma} 

\begin{lemma} \label{H0} If $G$ is a definably compact group and $H<G$ is a definable
subgroup, then $G^{00} \cap H < H^0$. \end{lemma}

\begin{proof} By definable choice (see \cite{vdDries98}) there is $D\subset G$
such that $G = \bigsqcup_{d\in D} Hd$ (disjoint union) $= HD$. Let
$h\in G^{00} \cap H$ and suppose for a contradiction that $h\nin
H^0$. Then $hH^0 \cap H^0 = \emptyset$. By the choice of $D$ it then
follows that $hH^0D \cap H^0D = \emptyset$. Since $H^0$ has finite
index in $H$, $H^0D$ is generic in $G$.  Together with $hH^0D \cap
H^0D = \emptyset$ this implies that $h\nin \stab_{\nog}(H^0D)$. So by
Corollary \ref{stab} $h\nin G^{00}$, a contradiction.  \end{proof}

\begin{theorem} \label{subgroup} If $G$ is a definably compact group and $H$
is a definable subgroup, then $H^{00} = G^{00} \cap H$. \end{theorem}

\begin{proof} The special case when $H$ is abelian and definably connected was proved
in \cite[Proof of Thm. 8.1, p. 33]{HrushowskiPP05}. The argument is the following. 
By Lemma \ref{torsionfreeabelian} it
suffices to show that $G^{00}\cap H$ is torsion free.  Fix a positive
integer $n$. Let $H[n] = \{h\in H \mid h^n = e\} < H$. By definable
choice there is a definable $X \subset H$ such that $H =
\bigsqcup_{h\in H[n]} hX$ and a definable $D \subset G$ such that $G =
\bigsqcup_{d \in D} Hd$. Hence $G = \bigsqcup_{h\in H[n]} hXD$. We
have $G^{00} \subset \stab_{\nog}(XD)$ and since the latter set
intersects $H[n]$ only in the identity element, $G^{00} \cap H[n] =
\{e\}$. Hence $G^{00} \cap H$ has no $n$-torsion.

The general case can be reduced to the special case as follows. By
Lemmas \ref{G0} and \ref{H0} we can assume that $H$ is definably
connected. Indeed from the definably connected case and the
Lemmas we deduce $G^{00} \cap H = G^{00} \cap H^0 = (H^0)^{00} =
H^{00}$. 

It remains to prove the case when $H$ is definably connected. By
\cite{Edmundo05} or \cite[Thm.6.12]{Berarducci05}, $H$ is the union of
its maximal definably connected abelian definable subgroups $T < H$
(which are moreover conjugates of each other). By the special case
$G^{00} \cap T = H^{00} \cap T = T^{00}$. Since the union of all the
$T$'s is $H$ we get $G^{00} \cap H = H^{00} = \bigcup_T T^{00}$. \end{proof}

\begin{remark} \label{tori} The proof also shows that if $G$ is a
definably compact definably connected group, then \[G^{00} = \bigcup_T
T^{00}\] where $T$ ranges over the maximal definably connected abelian
subgroups of $G$. Since each $T^{00}$ is divisible and torsion free,
we conclude that  $G^{00}$ is divisible and torsion free. \end{remark}

We can now eliminate the commutativity assumption in Lemma
\ref{torsionfreeabelian}.

\begin{theorem} Let $G$ be a definably compact definably connected group.  Then
$G^{00}$ is the unique type-definable subgroup of bounded index which
is torsion free. \end{theorem}

\begin{proof} By Remark \ref{tori} $G^{00} = \bigcup_T T^{00}$ where $T$ ranges
over the maximal definably connected abelian definable subgroups of
$G$. Each $T^{00}$ is torsion free. Hence so is their union
$G^{00}$. Conversely suppose $H$ is a torsion free type-definable
subgroup of $G$ of bounded index. If $T$ is a maximal definably
connected abelian definable subgroups of $G$, then $H\cap T$ has
bounded index in $T$ (as $[T:H\cap T] \leq [G:H]$) and is torsion
free. Therefore by Lemma \ref{torsionfreeabelian} $H\cap T = T^{00}$.
Taking the union over all $T$'s we obtain $H = G^{00}$. \end{proof}

Dropping the assumption that $G$ is definably connected we still
obtain a nice characterization:

\begin{corollary} \label{cabeliantf} Let $G$ be a definably compact group. 
Then $G^{00}$ is the unique type-definable
subgroup of bounded index which is divisible and torsion free.  \end{corollary}

\begin{proof} Let $H$ be a type-definable subgroup of $G$ of bounded index which
is divisible and torsion free. Since divisibility is preserved under
quotients, the image of $H$ under $\pi \colon G\to G/G^{0}$ is
divisible. But a divisible finite group must be trivial. So $H<G^0$
and by the previous theorem $H = (G^0)^{00} = G^{00}$. \end{proof} 

\begin{corollary} \label{Gprod}  Let $G$ be a definably compact group. 
Then $(G\times G)^{00} = G^{00} \times G^{00}$. \end{corollary}

This last corollary has the following consequence concerning the
measure $\mu_p$ induced by a generic type (see equation (2) in section 3).

\begin{proposition} Let $p$ be a generic type on $G$. There is a generic type $q$ on
 $G\times G$ such that for each pair of definable sets
 $A,B\subset G$ we have $\mu_{q}(A\times B) =
 \mu_p(A)\mu_p(B)$.  \end{proposition}

\begin{proof} Choose $q$ so that it contains all the sets whose complement is
non-generic and all the sets of the form $A\times B$ with $A,B\in
p$. This is possible since this collection of sets has the finite
intersection property.  Notice that $q$ is generic and for every
$A,B\in \Def(G)$ we have $A\times B\in q$ iff $A\in p$ and $B\in
p$. In fact one direction is by definition of $q$, and the other
follows easily from the fact that if a set is not in a type its
complement does. Now by definition $\mu_q(A\times B)$ is the Haar
measure of the set $\{(g,h)(G\times G)^{00} \mid (g,h)(A\times B)\in
q\}$. By the choice of $q$ and Corollary \ref{Gprod} this set is the
cartesian product of $\{gG^{00} \mid gA\in p\}$ and $\{hG^{00} \mid
hB\in p\}$. Now use the fact that the Haar measure on $G\times G$ is
the product measure of the Haar measure on $G$ \cite[(3) p. 263]{Halmos74}.  
\end{proof}

\section{Functorial properties}

Let $\varphi \colon G \to H$ be a definable homomorphism of definably
compact groups. Then $\varphi$ induces a morphism $F (\varphi) \colon
G/G^{00}\to H/H^{00}$ sending $gG^{00}$ to $\varphi(g)H^{00}$. This is
well defined since $\varphi (G^{00}) \subset H^{00}$. (Proof:
$\varphi^{-1}(H^{00})$ is a type-definable subgroup of $G$ of bounded
index, so it contains $G^{00}$.) It follows easily from the
definitions that $F(\varphi)$ is continuous with respect to the logic
topologies. Clearly $F$ preserves compositions and identity maps. So
we have:

\begin{proposition} $F \colon G\mapsto G/G^{00}$ is a functor from definably
compact groups and definable homomorphisms to compact Lie groups and
continuous homomorphisms. \end{proposition}

\begin{theorem} The functor $G\mapsto G/G^{00}$ is exact. \end{theorem}

\begin{proof} Consider the following commutative diagram.

\begin{center}
$\begin{array}{lllllllll}
 & & 0 & & 0 & & 0 & & \\
  & & \downarrow & & \downarrow & & \downarrow & & \\
0 & \longrightarrow & H^{00} & \longrightarrow & G^{00} 
& \longrightarrow & B^{00} & \longrightarrow & 0 \\
  & & \downarrow & & \downarrow & & \downarrow & & \\
0 & \longrightarrow & H & \stackrel{i}{\longrightarrow} & G & \stackrel{\pi} 
{\longrightarrow} & B & \longrightarrow & 0 \\
 & & \downarrow & & \downarrow & & \downarrow & & \\
0 & \longrightarrow & H/H^{00} & \longrightarrow & G/G^{00} 
& \longrightarrow & B/B^{00} & \longrightarrow & 0 \\
  & & \downarrow & & \downarrow & & \downarrow & & \\
 & & 0 & & 0 & & 0 & &
 \end{array}$
\end{center}

Suppose that the middle row is an exact sequence of definably compact
groups and definable homomorphisms. This means that $i$ is injective,
$\pi$ is surjective and the image of $i$ is the kernel of $\pi$. So
$i(H) \lhd G$ and $B \simeq G/i(H)$. 

We must prove that the bottom row is an exact sequence of Lie
groups. By a routine diagram chasing argument (``nine Lemma'') it
suffices to show that the top row is exact. We can assume that $H\lhd
G$ and $B = G/H$. In this case proving the exactness of the top row
amounts to show that $H^{00} = G^{00} \cap H$ and $B^{00} =
\pi(G^{00})$.  The first equality is Theorem \ref{subgroup}. For the
second note that $\pi(G^{00})$ is a type definable subgroup of $B$,
and since $\pi$ is onto it has bounded index in $B$. So $\pi(G^{00})$
contains $B^{00}$. For the opposite inclusion note that
$\pi^{-1}(B^{00})$ is a type-definable subgroup of bounded index of
$G$, so it contains $G^{00}$.  \end{proof}

\section{Spectral spaces}

\begin{definition} (\cite{Pillay88}) Let $M_0 \prec M$ be a small model.  Given an
$M_0$-definable set $X$, let $\widetilde X = {\widetilde X}(M_0)$ be
the set of types $S_X(M_0)$ (which can be identified with the
ultrafilters of $M_0$-definable subsets of $X$) with the following
{\bf spectral topology}: as a basis of open sets we take the sets of
the form $\widetilde U$ with $U$ a $M_0$-definable open subset of
$X$. Note that since $M_0$ is small and $M$ is saturated, every type
$p\in \widetilde X$ is realized in $M$ (i.e. there is $x\in M^n$ with
$x \in \bigcap_{X\in p} X$). \end{definition}

The spectral topology is coarser than the usual {\bf Stone topology}
(or {\bf contructible topology}) where one takes as a basis of open
sets all the sets of the form $\widetilde U$ with $U\subset X$ an
arbitrary $M_0$-definable set, not necessarily open.  The spectral
topology generalizes the real spectrum of semialgebraic geometry. More
precisely, if $M_0$ is a real closed field, then the set of $n$-types
over $M_0$ with the spectral topology coincides with the real spectrum
of $M_0[X_1, \ldots, X_n]$ (see \cite{CosteR82,Pillay88}).

The spectral topology is in general not Hausdorff, but it is
quasi-compact (every open covering has a finite subcovering) since it
is coarser than the constructible topology (which is compact). In fact
the spectral topology makes $\widetilde X$ into a {\bf spectral
space}, namely $\widetilde X$ has a basis of quasi-compact open sets
stable under finite intersections and every irreducible closed set is
the closure of a unique point (for the verification of these
properties in the o-minimal case see \cite[Lemma
1.1]{Pillay88}). Since we are working over an expansion of a field
every definable set $X$ is definably normal, namely every pair of
disjoint definable closed sets can be separated by disjoint definable
open sets (one can use the $M$-valued metric $|x-y|$). It then
follows:

\begin{lemma} (\cite[Thm. 2.12]{EdmundoJP05}) $\widetilde X$ is {\bf normal}, namely every
pair of disjoint closed sets can be separated by disjoint open sets.\end{lemma} 

\begin{proof} Every closed set of $\widetilde X$ is of the form $\bigcap_{i\in
 I} \widetilde A_i$ where each $A_i$ is an $M_0$-definable closed subset of
 $X$.  Consider two closed disjoint subsets $\bigcap_{i\in I}
 \widetilde {A_i}$ and $\bigcap_{j\in J} \widetilde {B_j}$ of
 $\widetilde X$, where $A_i,B_j$ are closed $M_0$-definable subsets of
 $X$. By the quasi-compactness of $\widetilde X$ there are finite
 subsets $I_0$ of $I$ and $J_0$ of $J$ with $\bigcap_{i\in I_0}
 \widetilde {A_i} \cap \bigcap_{j\in J_0} \widetilde {B_j} =
 \emptyset$. Since $A \mapsto \widetilde A$ is a boolean algebra
 isomorphism onto its image, we obtain $\bigcap_{i\in I_0} {A_i} \cap
 \bigcap_{j\in J_0} {B_j} = \emptyset$. Since $X$ is definably normal
 there are disjoint definable open subsets $U,V$ of $X$ with
 $\bigcap_{i\in I_0} {A_i} \subset U$ and $\bigcap_{j\in J_0} {B_j}
 \subset V$. But then $\widetilde U$ and $\widetilde V$ are disjoint
 open subsets of $\widetilde X$ which separate $\bigcap_{i\in I}
 \widetilde {A_i}$ and $\bigcap_{j\in J} \widetilde {B_j}$. \end{proof}

\begin{definition} (see \cite[Example 2.2]{HrushowskiPP05}) \label{M0} Let $X$ be a
definable set and let $E$ be a bounded type-definable equivalence
relation. In this situation we always assume that $M_0 \prec M$ is a
small model such that $X$ is defined over $M_0$, $E$ is type-defined
over $M_0$ and each equivalence class of $X/E$ contains an element
from $M_0$ (so $X/E \simeq X(M_0)/E(M_0)$). Let $\widetilde X =
S_X(M_0)$ with the spectral topology. There is a natural surjective map \[\Psi
\colon \widetilde X \to X/E\] defined as follows. Let $b\in X$ realize
the type $p$.  Then $\Psi (p) = \pi(b)$, where $\pi \colon X \to X/E$
is the natural projection. \end{definition}

Note that $\Psi$ is well defined since if $a,b$ have the same type
over $M_0$, then $\pi(a) = \pi(b)$. (If not by the choice of $M_0$
there is $a' \in X(M_0)$ with $E(a,a')$ and $\lnot E(b,a')$, hence
$E(x,a')$ is contained in the type of $a$ but not in the type of $b$.)
Moreover $\Psi$ is continuous with respect to the constructible
topology on $\widetilde X$ and the logic topology on
$X/E$ \cite[Example 2.2]{HrushowskiPP05}.

\begin{theorem} \label{spec1}
Assume $\pi \colon X \to X/E$ is continuous. Then $\Psi \colon
\widetilde X \to X/E$ is continuous with respect to the spectral
topology on $\widetilde X$ and the logic topology on $X/E$. \end{theorem}

\begin{proof}
Let $\pi \colon X \to X/E$ be the projection and let $Z\subset
X/E$ be a closed subset. Then by Lemma \ref{repres} we can write
\[\pi^{-1}(Z) = \bigcap_{i\in I} X_i\] where $I$ is small and for
each $i\in I\;$ $X_i$ is an $M_0$-definable set. (Requiring that $I$
is small is actually redundant since the family of all the
$M_0$-definable sets is small.)

\begin{claim}\label{cspec1}
\[\Psi^{-1}(Z) = \bigcap_i \widetilde {X_i}.\]
\end{claim}

In fact consider a type $p\in \widetilde X$ and let $b \in X$
realize $p$. Then $p\in \bigcap_i \widetilde {X_i}$ iff $b \in
\bigcap_i X_i$ iff $\pi(b) \in Z$ iff $p \in \Psi^{-1}(Z)$,
where the last equivalence follows from $\Psi(p) = \pi(b)$ (by
definition of $\Psi$).

Having proved the claim, to finish the proof of the main result we use
the assumption that $\pi$ is continuous. This ensures that
$\pi^{-1}(Z) = \bigcap_{i\in I} X_i$ is closed. We can further assume
that the family $\{X_i\mid i\in I\}$ is downward directed as otherwise
we consider its closure under finite intersections.  Since $I$ is
small and $M$ is saturated, we can easily conclude:

\begin{claim}\label{c2spec1}
\[\ov{\bigcap_i X_i} = \bigcap_i \ov{X_i}\]
\end{claim}

In fact it suffices to observe that if a definable open set $U$ is
disjoint from $\bigcap_i X_i$, then by saturation $U$ is disjoint from a
finite sub-intersection. By the claim we can assume that in the representation
$\pi^{-1}(Z) = \bigcap_i X_i$ each $X_i$ is an $M_0$-definable closed
subset of $X$. But then $\widetilde {X_i}$ is closed in the spectral
topology, and therefore so is $\Psi^{-1}(Z) = \bigcap_i \widetilde
{X_i}$. Thus $\Psi$ is continuous. \end{proof}

\begin{lemma} \label{closed} Assume $\pi \colon X \to X/E$ is continuous. Then
$\Psi \colon \widetilde X \to X/E$ is a closed map. \end{lemma}

\begin{proof} Let $Z\subset \widetilde X$ be closed. Then $Z$ is quasi-compact,
and therefore (since $\Psi$ is continuous) $\Psi(Z)\subset X/E$ is
quasi-compact. But $X/E$ is compact (Hausdorff), so the quasi-compact
subsets of $X/E$ are actually compact, hence closed. \end{proof}

\begin{corollary} \label{spec2}
Assume $\pi \colon X \to X/E$ is continuous. Then $\Psi \colon
\widetilde X \to X/E$ is a quotient map, namely a subset $Z\subset
X/E$ is open if and only if $\Psi^{-1}(Z)$ is open. \end{corollary}

\begin{proof} It suffices to show that a subset $Z\subset X/E$ is closed iff
$\Psi^{-1}(Z)$ is closed. One implication follows from the continuity
of $\Psi$ and the other from the fact that $\Psi$ is a closed
surjective map. \end{proof}

\begin{lemma} (\cite[Lemma 3.2]{Pillay04}) \label{cont} Let $G$ be a definable
group, and let $H<G$ be a type-definable subgroup of bounded
index. Then $H$ is open in $G$. So the natural map $\pi\colon G \to
G/H$ is continuous (indeed the preimage of any subset is open). \end{lemma}

By Lemma \ref{cont} and Corollary \ref{spec2} we have:

\begin{corollary} \label{quotient} $\Psi \colon \widetilde G \to G/H$ is a
continuous closed surjective map, hence a quotient map. \end{corollary}

By the corollary $G/H$ is homeomorphic to $\widetilde G/\ker
(\Psi)$ where $\ker (\Psi)$ is the equivalence relation $\{(p,q)
\mid \Psi(p) = \Psi(q)\}$ and $\widetilde G/\ker (\Psi)$ has the
quotient topology. Thus going to the spectral space we have
managed to understand $G/H$ in terms of the quotient topology.

\section{Cohomology}

Let $X\subset M^n$ be a definable closed and bounded set. By the
triangulation theorem $X$ is definably homeomorphic to the geometrical
realization $|K|$ (over $M$) of a finite simplicial complex $K$. There
are various ways one can define the cohomology group $H^n(X)$ (with
coefficients in $\Z$) but it is clear that for any reasonable choice
$H^n(X)$ should coincide up to isomorphism with the classical groups
$H^n(K)$. One possibility is to use sheaf cohomology. Namely one
defines $H^n(X)$ as $H^n(\widetilde X, Sh(\Z))$ where $Sh(\Z)$ is the
sheaf generated by the constant sheaf $\Z$ on $\widetilde X$ (see
\cite{Delfs80,DelfsK82,Delfs85,Delfs91,Knebusch92} in the
semialgebraic case and \cite{EdmundoJP05} in the o-minimal
case). Since $\widetilde X$ is a normal spectral space $H^n(\widetilde
X, Sh(\Z))$ is naturally isomorphic to the \v{C}ech cohomology group
$\check{H}^n(X;\Z)$ (see \cite[Prop. 5]{CarralC83}). For simplicity of
notation we use cohomology with coefficients in $\Z$, but all our
considerations apply to any other coefficient group.

\begin{remark} Equivalently one can work with sheaves directly on $X$ rather than
on $\widetilde X$ but then one has to consider $X$ not as topological
space but as a site in the sense of Grothendieck, which essentially
means that the only admissible open covers are the covers by finite
families of definable open sets. In this way the sheaves on $X$ are
naturally identified with those on $\widetilde X$
\cite[\S 1.3]{CarralC83}. \end{remark}

\begin{remark} If $X$ is a definably connected set, the spectral space
$\widetilde X$ is connected but in general it is not path
connected (not even locally). So it would not be satisfactory to
use singular cohomology for these spaces. \end{remark}

\begin{remark} 
It follows from the above discussion and the results in \cite{EdmundoJP05} that
\[\check{H}^n(\widetilde X; \Z) \simeq H^n_{\df}(X;\Z)\]
where $H^n_{\df}$ is the definably singular cohomology studied in
\cite{EdmundoO04} and based on the definable singular homology of
\cite{Woerheide96}.  
\end{remark}

\begin{definition} By Corollary \ref{quotient} there is a continuous closed
surjective map $\Psi\colon \widetilde G \to G/G^{00}$. Let $\Psi^*$ be
the induced homomorphism in \v{C}ech cohomology:
\[\Psi^* \colon \check{H}^n(G/G^{00};\Z) \to \check{H}^n(\widetilde
G;\Z)\] \end{definition}

Notice that $G/G^{00}$ is a very nice space (a compact Lie group),
so for this space \v{C}ech cohomology coincides with singular
cohomology. 

\begin{conjecture} \label{cohc}
$\Psi^*$ is an isomorphism. \end{conjecture}

\begin{remark} In the abelian case, by a result of Edmundo and Otero in
\cite{EdmundoO04}, we have \[H^n(G/G^{00};\Z) \simeq H^n_{\df}(G;\Z)\]
where $H^n$ is singular cohomology and $H^n_{\df}$ is definable
singular cohomology. By the above discussion, we can equivalently use
the \v{C}ech cohomology groups, thus obtaining
$\check{H}^n(G/G^{00};\Z) \simeq \check{H}^n(\widetilde G;\Z)$. This
however does not settle Conjecture \ref{cohc} even in the abelian
case. Indeed the proof in \cite{EdmundoO04} only tells us that there
is an isomorphism, but does not tell us which function does the
job. Notice in particular that there is no obvious way in which the
map $\pi\colon G \to G/G^{00}$ can induce an homomorphism
$H^n(G/G^{00};\Z) \to H^n_{\df}(G;\Z)$. \end{remark}

\begin{lemma} (\cite[Remark 2.17]{EdmundoJP05}) \label{ftilde} Let $f\colon X\to
Y$ be an $M_0$-definable definable function. Then $f$ induces a
function $\widetilde f \colon \widetilde X \to \widetilde Y$ by $f(p)
= \{Z \mid f^{-1}(Z) \in p\}$ where $Z$ ranges over the
$M_0$-definable subsets of $Y$.  We have $\widetilde{f(X)} =
\widetilde f \widetilde X$. Moreover if $f$ is continuous, $\widetilde
f$ is continuous. So if $f$ is an homeomorphism, then $\widetilde f$
is an homeomorphism. \end{lemma}

\begin{definition} Let $\widetilde {G^{00}}\subset \widetilde G$ be the set of types
which are realized by elements in $G^{00}$. Equivalently, if $G^{00} =
\bigcap_i X_i$ with $X_i$ $M_0$-definable, $\widetilde {G^{00}}= \bigcap_i
\widetilde {X_i}$ (this does not depend on the choice of the
representation).  \end{definition}

\begin{lemma} \label{homeo} Each fiber $\Psi^{-1}(hG^{00})$ of $\Psi\colon \widetilde G \to
G/G^{00}$ is homeomorphic to $\widetilde{G^{00}}$. \end{lemma}

\begin{proof} Let $G^{00} = \bigcap_i X_i$ with $X_i$ $M_0$-definable.  Let $p\in
\widetilde G$ and let $g\in G$ realize $p$. So by definition $\Psi(p)
= gG^{00}$. By the choice of $M_0$ (see Definition \ref{M0}) each
coset of $G/G^{00}$ is of the form $hG^{00}$ with $h$ in $M_0$. We
have the following chain of equivalences: $p \in \Psi^{-1}(hG^{00})$
iff $gG^{00} = hG^{00}$ iff $g\in h(\bigcap_i X_i) =\bigcap_i h X_i$
iff $p\in \bigcap_i \widetilde {h X_i}$. So $\Psi^{-1}(hG^{00}) =
\bigcap_i \widetilde {h X_i}$. Let $h\colon G \to G$ be left
multiplication by $h$ and let $\widetilde h \colon \widetilde G\to
\widetilde G$ be the induced map as in Lemma
\ref{ftilde}. Then $\Psi^{-1}(hG^{00}) = \bigcap_i \widetilde h
\widetilde {X_i} = \widetilde h \bigcap_i \widetilde {X_i} =
\widetilde h \widetilde {G^{00}}$. So $\Psi^{-1}(hG^{00})$ is
homeomorphic via $\widetilde h$ to $\widetilde {G^{00}}$. \end{proof}

Our next goal is to investigate the cohomology of $\widetilde {G^{00}}$. 

\begin{lemma} \label{qcomp} $\widetilde {G^{00}}$ is a closed subset of
$\widetilde G$, so it is quasi-compact. \end{lemma}

\begin{proof} $G/G^{00}$ is a compact Hausdorff space, so in particular its
points are closed.  By Lemma \ref{homeo} $\widetilde {G^{00}}$ is the
preimage of a point of $G/G^{00}$ under a continuous map, so it is
closed. \end{proof}

The fact that $\widetilde {G^{00}}$ is closed also follows from the
following lemma.

\begin{lemma} \label{Cl} Suppose $G^{00} = \bigcap_i X_i$ where $\{X_i\mid i\in
I\}$ is a downward directed small family of definable sets. Then
$G^{00} = \bigcap_i Cl(X_i)$. \end{lemma}

\begin{proof} By Lemma \ref{ctble} $G^{00} = \bigcap_{n\in \NN} Y_i$ with
$\Cl(Y_{n+1}) \subset \Int(Y_n)$. Given $n\in \NN$ by saturation a
finite sub-intersection of the $X_i$'s is contained in $Y_n$, and since
the family is directed there is a single $X_{i_n}$ contained in
$Y_n$. Now $\bigcap_i Cl(X_i) \subset \bigcap_n Cl(X_{i_n}) \subset
\bigcap_n \Cl(Y_n) = G^{00}$. \end{proof}

\v{C}ech cohomology enjoys the following continuity property. If $Y$
is a compact subset of a normal space $X$, then the \v{C}ech
cohomology of $Y$ is the inductive limit of the \v{C}ech cohomologies
of its open neighbourhoods in $X$. A readable proof of this fact can
be found in \cite[Lemma 73.3]{Munkres84}. Exactly the same argument
shows:

\begin{lemma} \label{inductive} 
If $Y$ is a closed subset of a normal quasi-compact space $X$,
then the \v{C}ech cohomology of $Y$ is the inductive limit of the
\v{C}ech cohomologies of its open neighbourhoods in $X$. \end{lemma}

A more general statement along the same lines, but using
sheaf cohomology, is the following:

\begin{lemma} (\cite[Thm. 3.1]{Delfs85}) \label{delfs} Let $\cal F$ be an
arbitrary sheaf on the normal spectral space $\widetilde X$ and let
$Y$ be a quasi-compact subset of $\widetilde X$. Then the canonical
homomorphism
\[\varliminf_{\substack{Y\subset U \\ U \mbox{ open in } 
\widetilde X}} H^q(U;{\cal F})\to H^q(Y;{\cal F}|Y)\] is an
isomorphism for every $q \geq 0$. \end{lemma}

\begin{remark}
 \label{taut} In the terminology of \cite[Def. 10.5]{Bredon97} this implies that
$Y$ is a {\bf taut} subspace of $\widetilde X$ (see
\cite[Thm. 10.6]{Bredon97}).\end{remark}

In \cite[Thm. 3.1]{Delfs85} the above result is stated under the
hypothesis that $\widetilde X$ is the the real spectrum of a ring, but
as observed in \cite{EdmundoJP05} the argument actually works for
every normal spectral space.

We plan to apply Lemma \ref{inductive} to the closed
(hence quasi-compact) subset $\widetilde {G^{00}}$ of the normal
spectral space $\widetilde G$. In the computation of the inductive
limits, we can restrict ourselves to a cofinal sequence of open
neighbourhoods of $\widetilde {G^{00}}$. To produce such cofinal
sequences we use the following result.

\begin{lemma} \label{cofinal} Suppose $G^{00} = \bigcap_i X_i$ where 
$\{X_i\mid i\in I\}$ is a downward directed family of $M_0$-definable
sets.  Then each open neighbourhood $U$ of $\widetilde{G^{00}}$ in
$\widetilde G$ contains one of the sets $\widetilde {X_i}$. So if the
$X_i$'s are open, then
\[\check{H}^n(\widetilde {G^{00}};\Z) =
\varliminf_{\; i} \check{H}^n(\widetilde {X_i};\Z).\] 
\end{lemma}

\begin{proof} By Lemma \ref{Cl} $G^{00} = \bigcap_i Cl(X_i)$. So $\bigcap_i
\widetilde{Cl(X_i)} \subset U$. Since $\widetilde G$ is quasi-compact,
there is a finite sub-intersection of the closed sets
$\widetilde{Cl(X_i)}$ which is contained in the open set $U$. Since
the family $\{X_i\}$ is directed, a single $\widetilde{Cl(X_i)}$ is
contained in $U$. So a fortiori $\widetilde {X_i}$ is contained in
$U$. \end{proof}

In the next section we consider the conjecture that $\widetilde
{G^{00}}$ is acyclic in \v{C}ech cohomology. Notice that Lemma
\ref{cofinal} leaves open even the question whether the \v{C}ech
cohomology of $\widetilde {G^{00}}$ is finitely generated.

\section{Contractibility}

Let $G$ be a definably compact group. By Lemma \ref{ctble} $G^{00}$ is
the intersection of a decreasing sequence of definable open subsets of
$G$. 

\begin{conjecture} \label{contr}
$G^{00}$ is the intersection $\bigcap_{i\in \NN} X_i$ of a
decreasing sequence $X_0\supset X_1\supset X_2 \supset \ldots$ of
 $M_0$-definable definably contractible open sets. \end{conjecture}

\begin{proposition} \label{SO} 
The conjecture holds for the group $G = SO(n,M)$ of Example
\ref{Example}. \end{proposition}

\begin{proof} The idea is to use transfer from $\R$ to $M$. Let $T_e$ be the
tangent plane of $SO(n,\R)$ at the neutral element $e$ of the
group. Every sufficiently small ball around the origin in $T_e$ is
semialgebraically homeomorphic to a neighbourhood of $e$ in $SO(n,\R)$
via a projection $p$ along the orthogonal complement of $T_e$ in the
ambient space $\R^{n^2}$. So we can fix a positive integer $k$ such
that $p$ maps the ball of radius $1/k$ in the tangent plane
homeomorphically onto a closed neighbourhood $X_k$ of $e$ in
$SO(n,\R)$. Moreover there are balls $B_k\subset C_k \subset \R^{n^2}$
centered at $e$ and of rational radius $b_k$ and $c_k$ respectively,
such that $b_k$ and $c_k$ tend to $0$ with $k\to \infty$ and $B_k \cap
SO(n,\R) \subset X_k \subset C_k \cap SO(n,\R)$.  Once $k, b_k, c_k$
have been fixed, these statements are semialgebraic, so they transfer
to $SO(n,M)$ by the Tarski-Seidenberg principle. Since $SO(n,M)^{00}$
is the kernel of the standard part map $\sta \colon SO(n,M) \to
SO(n,\R)$ it follows that $\bigcap_k X_k(M) = SO(n,M)^{00}$. Moreover
$X_k(M)$ is semialgebraically contractible since it is
semialgebraically homeomorphic to a ball in the tangent plane. \end{proof}

With the same argument one can prove the conjecture for any
definable group with ``good reduction'' in the sense of
\cite{HrushowskiPP05}. We recall the definition: 

\begin{definition} $G$ has {\bf good reduction} if
it is definably isomorphic in $M$ to a group $G_1$ which can be
defined over $\R$ in the following sense: there is a sub-language
$L_0$ of the language of $M$ which contains $+,\cdot$ and there is
an elementary substructure $M_0$ of $M|L_0$ of the form $M_0 =
(\R,+,\cdot, \ldots)$ such that $G$ is defined over $M_0$. \end{definition}

The arguments in Proposition \ref{SO} establish the following:

\begin{theorem} Suppose $G$ is a definably compact group with good reduction. Then
the contractibility conjecture \ref{contr} holds for $G$.  \end{theorem}

In particular the contractibility conjecture \ref{contr} holds for
every definably simple group $G$, since such groups have good
reductions \cite[Thm. 5.1]{PeterzilPS02}. Another favorable case is the following. 

\begin{theorem} If $G$ is a definably compact group with $\dim(G) = 1$, then 
Conjecture \ref{contr} holds for $G$. \end{theorem}

\begin{proof} By Lemma \ref{G0} we can assume that $G$ is definably connected.
By \cite{Pillay88b} $G$ is abelian and it is a definable manifold
without boundary. By \cite{Razenj91} (see also
\cite[Prop. 3.5]{Pillay04}) for each $x\in G$ there is a definable
linear ordering $<_x$ on $G\setminus \{x\}$ such that the topology of
$G\setminus \{x\}$ coincides with the topology generated by the open
intervals of $<_x$ (recall that we always assume that $G$ is embedded
in some $M^n$ with the induced topology). Moreover $G$ has a unique
element $x$ of order $2$ and by \cite[Prop. 3.5]{Pillay04} $G^{00}$ is
the intersection of all the $<_x$-intervals $[-a,a] \subset G\setminus
\{x\}$ where $a\in G$ ranges over the torsion elements of $G$ (with
the exclusion of $x$). By \cite[Thm. 6.5]{BerarducciO01} a
one-dimensional definably connected definably compact manifold with
non-empty boundary is definably homeomorphic to $[0,1]\subset M$. A
similar argument shows that $G\setminus \{x\}$ is definably
homeomorphic to $(0,1)\subset M$. So the $<_x$-intervals are definably
homeomorphic to intervals of $(0,1)$ and therefore are definably
contractible. \end{proof}

Let us now study some consequences of Conjecture \ref{contr}. We need: 

\begin{lemma} \label{homot} Suppose $X$ is an $M_0$-definable definably contractible
set. Then $\widetilde X$ is acyclic in \v{C}ech cohomology, namely
$\widetilde X$ is connected and $\check{H}^n(\widetilde {X};\Z)= 0$
for all $n>0$. \end{lemma} 
\begin{proof} By the verification of the homotopy axiom in
\cite[Thm. 4.1]{Delfs85} in the semialgebraic case and in
\cite{EdmundoJP05} in the o-minimal case. \end{proof}

\begin{theorem} \label{zero} Let $G$ satisfy the contractibility conjecture
\ref{contr}. Then $\widetilde {G^{00}}$ is acyclic in \v{C}ech
cohomology, namely $\widetilde {G^{00}}$ is connected and
$\check{H}^n(\widetilde {G^{00}};\Z) =0$ for all $n>0$. \end{theorem}

\begin{proof} By Lemmas \ref{homot}, \ref{cofinal} and \ref{homeo}. \end{proof}

\begin{corollary} \label{vietoris} Suppose $\widetilde {G^{00}}$ is acyclic in
\v{C}ech cohomology. Then \[\Psi^* \colon \check{H}^n(G/G^{00};\Z) \simeq
\check{H}^n(\widetilde {G};\Z).\] \end{corollary}

\begin{proof} The idea is to use the Vietoris-Begle mapping theorem. A simple
version of the theorem says that if $A$ is compact and $f\colon A\to
B$ is a surjective continuous map with acyclic fibers (with respect to
\v{C}ech or equivalently Alexander-Spanier cohomology), then $f$
induces an isomorphism in cohomology \cite{Spanier81}.  We cannot
directly apply this version of the theorem since $\widetilde G$ is
quasi-compact and normal but not Hausdorff. A possible way to handle
this problem is to work with the subspace ${\widetilde G}^{\max}$ of
the closed points of $\widetilde G$, which is compact Hausdorff and
has the same \v{C}ech cohomology. Alternatively we can use the
following more general version of the Vietoris-Begle mapping theorem
which works for sheaf cohomology and arbitrary topological spaces (we
give the version without family of supports):

\begin{lemma} (\cite[Ch. 2, Thm. 11.7]{Bredon97})\label{bredon} Let $f\colon
X\to Y$ be a closed continuous surjection, let $\cal B$ be a sheaf on
$Y$. Also assume that each $f^{-1}(y)$ is connected and taut in $X$
and that $H^p(f^{-1}(y); f^*{\cal B}_y) = 0$ for $p>0$ and all $y\in
Y$. Then
\[f^* \colon H^*(Y;{\cal B}) \to H^*(X;f^*{\cal B})\]
is an isomorphism. \end{lemma}

We can apply Lemma \ref{bredon} to the map $\Psi \colon \widetilde G
\to G/G^{00}$ and the constant sheaf ${\cal B} = Sh(\Z)$ over
$G/G^{00}$.  Then the inverse image $f^*{\cal B}$ is the constant
sheaf $Sh(\Z)$ over $\widetilde G$. We have already remarked that in
this case sheaf-cohomology coincide with \v{C}ech cohomology (over
$\Z$). The hypothesis of Lemma \ref{bredon} are verified by $\Psi$
thanks to Theorem \ref{zero} and Remark \ref{taut}. This proves
Corollary \ref{vietoris} \end{proof}

An attempt to prove the contractibility conjecture \ref{contr}, and
therefore the acyclicity of $\widetilde {G^{00}}$, could be based on
the following.

\begin{proposition} 
Suppose that for each $M_0$-definable set $X\subset G$ there is a
$M_0$-cell decomposition of $X$ such that, for each cell $C$ of the
decomposition, $CC^{-1}$ is $M_0$-definably contractible. Then
Conjecture \ref{contr} holds for $G$. \end{proposition}

\begin{proof} By Lemma \ref{ctble} $G^{00}$ is the intersection of a decreasing
sequence of $M_0$-definable open sets $O_i$ ($i\in \NN$).  Since $O_i
\supset G^{00}$, $O_i$ is generic.  Since $G^{00} = \bigcap_i O_i$ is
a group, by saturation $\bigcap_i O_iO_i^{-1} = (\bigcap_i
O_i)(\bigcap_i O_i)^{-1} = G^{00}$.  Built inductively cell
decompositions ${\cal D}_i$ of $O_i$ refining the previous and such
that for all $C\in {\cal D}_i$, $CC^{-1}$ is definably contractible.
At least one of the cells of $O_i$ is generic. Since each cell of
$O_i$ is a union of cells of $O_{i+1}$, each generic cell of $O_i$
contains a generic cell of $O_{i+1}$. By K\"onig's lemma there is an
infinite sequence of generic cells $C_i \subset O_i$ with $C_{i+1}
\subset C_i$. Although $C_i$ does not need to contain $G^{00}$, by
Lemma \ref{CC} we have $G^{00} \subset C_i{C_i}^{-1}$. On the other
hand $\bigcap_i C_iC_i^{-1} \subset \bigcap_i O_iO_i^{-1} \subset
G^{00}$, so $G^{00} = \bigcap_i C_iC_i^{-1}$ and since by construction
$C_iC_i^{-1}$ is definably contractible we are done. \end{proof}

\subsection*{Acknowledgments}  
I thank Antongiulio Fornasiero and Fulvio Lazzeri for many
conversations on the topics of this paper. In particular countless discussions with 
Lazzeri helped me filtering various ideas and conjectures on the cohomological
aspects. Part of results of this paper were presented in a preliminary
form at the meeting ``Around o-minimality'', March 11-13, 2006, Leeds,
organized by Anand Pillay.

\end{document}